\documentclass[11pt]{amsart}
\usepackage{amsfonts}

\newtheorem{theorem}{Theorem}[section]
\newtheorem{lemma}[theorem]{Lemma}

\newtheorem{definition}[theorem]{Definition}
\newtheorem{corollary}[theorem]{Corollary}

\newcommand{\N}{\mbox{$\mathbb{N}$}}

\newcommand{\K}{\mbox{$\mathbb{K}$}}

\newcommand{\A}{\mbox{${\mathcal A}$}}
\newcommand{\B}{\mbox{${\mathcal B}$}}
\newcommand{\C}{\mbox{${\mathcal C}$}}
\newcommand{\D}{\mbox{${\mathcal D}$}}
\newcommand{\E}{\mbox{${\mathcal E}$}}
\newcommand{\F}{\mbox{${\mathcal F}$}}
\newcommand{\G}{\mbox{${\mathcal G}$}}
\newcommand{\Li}{\mbox{${\mathcal L}$}}
\newcommand{\M}{\mbox{${\mathcal M}$}}
 
\begin{document}
\title[Morita theorem for Operator algebras]
{A Morita theorem for Algebras of Operators on 
Hilbert Space}
 
\vspace{30 mm}
 
\author{David P. Blecher}
\date{September, 1998.}
\thanks{* Supported by a grant from the NSF}
\address{Department of Mathematics\\University of Houston\\Houston,
TX 77204-3476 }
\email{dblecher@@math.uh.edu}\maketitle
 
\maketitle
  
\vspace{10 mm}
 
\begin{abstract}
We show that two operator algebras 
are strongly Morita equivalent (in the sense of Blecher, Muhly 
and Paulsen) if and only if
their categories of operator modules 
are equivalent via completely
contractive functors.  Moreover, any such functor 
is completely isometrically
isomorphic to the Haagerup tensor product 
(= interior tensor product) with
a strong Morita equivalence bimodule.  
\end{abstract}
 
\pagebreak
\newpage
 
\setcounter{section}{0}

\section{Notation, background and statement of the theorem}
Around 1960, in pure algebra, arose the notion of 
Morita equivalence of rings.  Two rings $A$ and $B$ were 
defined to be Morita equivalent if
the two categories
$_{\A}MOD$ and $_{\B}MOD$ of modules
are equivalent.
The fundamental theorem from 
those early days of that subject 
(see \cite{F,Morita,Bass}) is that these categories
are equivalent if and only if there exists a
pair of bimodules $X$ and $Y$ such that
$X \otimes_B Y \cong A$ and $Y \otimes_A X \cong B$ 
as bimodules.  The theorem goes on to describe 
these so-called `equivalence bimodules'
 and how they arise, and the implications
for $\A$ and $\B$.   

In the early 70's M. Rieffel introduced and 
developed the notion of
{\em strong Morita equivalence of 
C$^*-$algebras} (see \cite{Ri2} for a 
good discussion and survey).  It has become a fundamental
tool in modern operator algebra and noncommutative geometry
(see \cite{Con} for example).  
Rieffel defined strong Morita equivalence in terms
of the existence of a certain type of bimodule, possessing
certain C$^*-$algebra valued positive definite inner products.
Until recently there was no description in terms of
a categorical equivalence.  Except for the absence
of such a theorem, the basic features from pure algebra
were shown to
carry over quite beautifully.  Of course one expects,
and obtains, stronger (functional analytic) variants
of these basic features.  As just one example:
 in pure algebra, one finds that $\B \cong p M_n(\A) p$,
for a projection $p$ in the $n \times n$ matrices
$M_n(\A)$.
The same thing is true in the case of 
unital strongly Morita equivalent C$^*-$algebras
$\A$ and $\B$, except that $p$ is an orthogonal 
projection, and the `$\cong$' means `as C$^*-$algebras', 
i.e. (isometrically)  $*$-isomorphically. 

In \cite{BlFTM} we 
showed that two C*-algebras
are strongly Morita equivalent if and only if
their categories of 
(left) operator modules (defined below)
are equivalent via completely
contractive functors.  Moreover, 
any such functor is completely isometrically
isomorphic to 
the Haagerup tensor product 
(= interior tensor product \cite{L2,Bna}) 
with 
an equivalence bimodule.

Here we generalize this result to possibly nonselfadjoint
operator algebras, 
that is, to general norm closed algebras of operators 
on Hilbert space\footnote{We note that
the work of Morita on purely algebraic 
equivalence,
and many related consequences,
was summarized and popularized 
by Bass as a collection of theorems
known as Morita I, II and III (see \cite{Bass,F}).
Most of the appropriate version of `Morita I' was
proved for C$^*-$algebras by Rieffel, and for
general operator algebras in \cite{BMP}.  Our main theorem 
here is a `Morita II' theorem for (possibly nonselfadjoint)
operator algebras.
The appropriate version of  
`Morita III' follows easily from what we have
done, as in pure algebra, and is omitted.}.
Thus we answer the main remaining theoretical
question from our study of Morita equivalence of possibly 
nonselfadjoint operator algebras, begun in \cite{BMP}.
The various ingredients of
our proof shows how the algebra and functional analytic 
structures, in particular, the geometry of the
associated Hilbert spaces, are intricately connected.
Some major tools, such as von Neumann's double commutant 
theorem, do not exist for nonselfadjoint operator algebras;
to overcome this
we use the theory of C$^*-$dilations 
of operator modules developed in \cite{Bdil}, to transfer 
the problem
to the C$^*-$algebra scenario, where we may more or less 
use our earlier proof of \cite{BlFTM}, and the 
lowersemicontinuity argument on the quasistate space
which we used there.  Of necessity some 
of our argument consists of instructions on
how to follow along and adapt
steps in the proof in \cite{BlFTM}.  
In order to not try the readers patience more
than needs be, we attempted to
keep these instructions minimal, yet sufficient.

Let us begin 
by establishing the common symbols and notations in this paper.  
We shall use operator spaces and completely bounded maps
quite extensively, and their connections to 
operator algebras, operator modules and C$^*-$modules.  
We refer the reader to 
\cite{BK,P,BMP,Bna,BlFTM,Bdil}
for missing background.  It is perhaps worth saying 
to the general reader
that it has been clear for some time that to understand
a general operator algebra $\A$ or operator module, 
it is necessary not only 
to take into account the norm, but also the natural 
norm on $M_n(\A)$.  That is one of the key perspectives
of operator space theory.  Hence we are not interested 
in bounded linear transformations, rather we look 
for the completely bounded, completely isometric,
or completely contractive maps - where the adjective
`completely' means that we are applying our maps to
matrices too.  This is explained at length in the 
references mentioned above.
The algebraic background needed 
may be found in any account of Morita theory
for rings, such as \cite{AF} or \cite{F}.  

We will use the symbols  
 ${\mathcal A}, {\mathcal B}$ for operator algebras.
We shall assume that our operator algebras
have contractive
approximate identities (c.a.i.'s).  
It is well known 
that every C$^*-$algebra is an operator algebra in this
sense.  
We write ${\mathcal C}$ and 
${\mathcal D}$ for the universal
or maximal C$^*-$algebras generated by
 ${\mathcal A}$
and ${\mathcal B}$
respectively, see \S 2 of \cite{Bdil}.
The symbol $r_v$ will always mean the 
`right multiplication by $v$'
operator, namely  $x \mapsto xv$, whose domain is usually 
the algebra $\A$ or $\C$.  
We will use the letters
 $H$ and $K$ for Hilbert spaces,  
$\zeta, \eta$ are
 typical elements in $H$ and
$K$ respectively, 
and $B(H)$ (resp. $B(H,K)$) is the
 space of bounded linear operators on 
$H$ (resp. from $H$ to $K$).  
 
Suppose that $\pi$ is a 
completely contractive representation of $\A$ on 
a Hilbert space $H$, and that $X$ is a closed subspace of 
$B(H)$ such that $\pi(\A) X \subset X$.  
Then $X$ is a left $\A$-module.  
We shall assume that the module action is nondegenerate 
(= essential). 
We say that such $X$, considered as
an abstract operator space and a left $\A$-module, is a
left {\em operator module} over $\A$.  
By considering $X$ as an abstract operator 
space and module, we may forget about the particular $H, \pi$ 
used.
An obvious modification of a
 theorem of Christensen, Effros 
and Sinclair \cite{CES} tells us that
the operator modules are (up to completely isometric
isomorphism) exactly the operator spaces $X$ which are 
(nondegenerate) left $\A$-modules, such
that the module action 
is a `completely contractive' bilinear map 
(that is $\Vert a x \Vert \leq \Vert a \Vert \Vert x \Vert$
for matrices $a$ and $x$ with entries in $\A$ and $X$ respectively)
or equivalently,
the module action linearizes to a 
complete contraction $\A \otimes_h X \rightarrow X$,
where $\otimes_h$ is the Haagerup 
tensor product).  Such an $X$ is referred to as
an abstract operator module. 
We will use the facts
that submodules and  quotient modules
of operator modules, are again 
operator modules.  
We write $_{\A}OMOD$ for the category of 
left $\A$-operator modules.  
The morphisms 
are $_{\A}CB(X,W)$, the {\em completely bounded} left 
$\A-$module maps.  If $X$ is also a right
$\B-$module, then $_{\A}CB(X,W)$ is a left $\B$-module
where $(bT)(x) = T(xb)$, or equivalently,
$bT = T r_b$.  We will write $_{\A}CB^{ess}(X,W)$
for the subset consisting of such maps $bT$, for $b \in B$.
If $X,W \in $ $_{\A}OMOD$ then $_{\A}CB(X,W)$ is an 
operator space \cite{ERbimod}.  In this paper, when 
$X, W$ are operator modules or bimodules, and when we say 
`$X \cong W$' , or `$X \cong W$ as operator modules', we will
mean that the implicit isomorphism is a completely isometric
module map.   

We will need the following important
principle from \S 3 of \cite{Bdil}
which we shall use several times here without
comment: an isometric surjective $\A$-module map between two
Banach $\C$-modules, is a $\C$-module map.  
This shows that the `forgetful functor' 
$
_{\C}OMOD \rightarrow $$ _{\A}OMOD$, 
embeds $_{\C}OMOD$ as a (non-full) subcategory of
$_{\A}OMOD$.
To an algebraist, it may be helpful to remark that
it is a reflective subcategory
in the sense of \cite{Cat}.  The C$^*$-dilation, or maximal
dilation, referred to earlier, is the left adjoint of this
forgetful functor; and it can be explicitly described as
the functor 
$\C \otimes_{h\A} -$.  Here $\otimes_{h\A}$ is the module 
Haagerup tensor product studied in \cite{BMP}.  We will 
repeatedly use the fact 
(3.11 in \cite{Bdil}) that the `obvious map'
$V \rightarrow \C \otimes_{h\A} V$, is completely 
isometric, thus $V$ is an $\A$-submodule of its
maximal dilation.

We now turn to the category  $_{\A}HMOD$ 
of Hilbert spaces $H$ 
which are left $\A-$modules via a nondegenerate
completely contractive 
representation of $\A$ on $H$.  If $\A$ is a C$^*-$algebra,
then this is the same as the category of nondegenerate
$*-$representations of $\A$ on Hilbert
space.   By the universal property of the maximal
generated C$^*-$algebra, $_{\A}HMOD = _{\C}HMOD$ as 
{\em objects}.
In \cite{BMP}
we showed how $_{\A}HMOD$ may be viewed as a subcategory of
$_{\A}OMOD$ (see the discussion 
at the end of Chapter 2, and after 
Proposition 3.8, there).
Briefly, if $H \in \; _{\A}HMOD$ then if $H$ is equipped with its
`Hilbert column' operator space structure $H^c$, then 
$H^c \in \; _{\A}OMOD$.
Conversely, if $V \in \; 
_{\A}OMOD$ is also a
 Hilbert
column space, then the associated representation $A \rightarrow
B(V)$ is completely contractive and nondegenerate.
It is well known that for a linear map $T : H \rightarrow K$
between Hilbert spaces, the
usual norm equals the completely bounded norm of $T$
as a map $H^c \rightarrow K^c$.
Thus we see that the assignment 
$H \mapsto H^c$  embeds $_{\A}HMOD$ as a full
subcategory of 
$_{\A}OMOD$.  In future, if a Hilbert space is referred to
as an operator space, it will be with respect to its
column operator space structure, unless specified to
the contrary.  

We are concerned with functors between categories of
operator modules.  Such functors   
$F : $ $_{\A}OMOD \rightarrow $ $_{\B}OMOD$ are  
assumed to be  linear on spaces of morphisms.  Thus 
$T \mapsto F(T)$, from the space $_{\A}CB(X,W)$ to 
$_{\B}CB(F(X),F(W))$, is linear, 
for all pairs of objects $X,W \in $ $_{\A}OMOD$.  
We say $F$ is  {\em completely contractive}, if
this map $T \mapsto F(T)$
is completely contractive, 
for all pairs of objects $X,W \in $ $_{\A}OMOD$.
We say two functors  
$F_1, F_2 : $ $_{\A}OMOD \rightarrow $ $_{\B}OMOD$ 
are (naturally) completely isometrically 
isomorphic, if they are naturally isomorphic in
the sense of category theory 
\cite{AF,Cat}, with the natural transformations being 
complete isometries.  In this case we write 
$F_1 \cong F_2$ {\em completely isometrically}.  
 
\begin{definition}
\label{defome} We say that two operator algebras 
$\A$ and $\B$ are (left) {\em operator Morita equivalent}
if there exist completely contractive  
functors $F : $ $_{\A}OMOD \rightarrow $ $_{\B}OMOD$
and $G : $ $_{\B}OMOD 
\rightarrow $ $_{\A}OMOD$, such that $F G \cong Id$ and 
$G F \cong Id$ completely isometrically.  
Such  $F$ and $G$ will be called {\em operator 
equivalence functors}.
\end{definition}    

There is an obvious adaption to 
`right operator Morita equivalence',
where we are concerned with right operator modules.  We remark
that for C$^*-$algebras it is easy to show that
left operator Morita equivalence 
implies right operator Morita equivalence,
but this seems much harder for nonselfadjoint operator algebras,
although we shall see that it is true.

In \cite{BMP} we generalized 
strong Morita equivalence of C$^*-$algebras 
to possibly nonselfadjoint operator algebras: 

\begin{definition}
\label{defsme}   Two operator algebras
$\A$ and $\B$ are {\em strongly Morita equivalent}
if there exists an $\A-\B-$operator bimodule $X$, and a
$\B-\A-$operator bimodule $Y$, such that 
$X \otimes_{h\B} Y \cong \A$ completely isometrically and
as $\A-\A$-bimodules, and such that 
$Y \otimes_{h\A} X \cong \B$ completely isometrically and
as $\B-\B-$bimodules.   We say that $X$ is 
an $\A-\B-${\em strong Morita equivalence
bimodule}.
\end{definition}

The above is not quite the definition given in 
\cite{BMP}, although
we remarked, without giving a proof, that it is an
equivalent definition.  
Essentially
it is the same proof of the corresponding result in pure algebra
(see \cite{Cohn} or \cite{F} 12.12.3 and 12.13).  In our 
scenario there is really
only one new point, 
that is the element $u$ described in these texts
is not
in $\A$, but in $CB_{\A}(\A,\A) = RM(\A)$, the right 
multipliers of $\A$.  However it
commutes with $\A$ in $RM(\A)$, so it falls in the center of the 
multiplier
algebra  $M(\A)$ of $\A$, and in fact it is a 
unitary there.  The rest of
the proof carries through quite obviously.

We can now state our main theorem.  

\begin{theorem}
\label{FT}  
Two operator algebras $\A$ and $\B$ with contractive approximate
identities
are strongly Morita equivalent if and only if they
are left operator Morita equivalent, and if and only if they
are right operator Morita equivalent.  Suppose that
 $F, G$ are the left operator equivalence 
functors, and set $Y = F(\A)$ and $X = G(\B)$.  Then $X$ is an
$\A-\B-$strong Morita equivalence bimodule,  and
$Y$ is a
$\B-\A-$strong Morita equivalence bimodule, which is
 unitarily equivalent 
to the dual operator module $\tilde{X}$ of $X$.  Moreover, 
$F(V) \cong Y \otimes_{h\A} V \cong $ $_{\A}{\mathbb K}(X,V)$  
completely isometrically 
isomorphically (as $\B-$operator modules), for all $V 
\in $ $_{\A}OMOD$.  
Thus $F \cong Y \otimes_{h\A} - \; \cong $ 
$_{\A}{\mathbb K}(X,-)$ completely isometrically.   
Similarly $G \cong X \otimes_{h\B} - \; \cong $ 
$_{\B}{\mathbb K}(Y,-)$ 
completely isometrically.   
Also $F$ and $G$ restrict to equivalences of the subcategory
$_{\A}HMOD$ with $_{\B}HMOD$, the subcategory $_{\C}HMOD$ with 
$_{\D}HMOD$, and
the subcategory $_{\C}OMOD$ with $_{\D}OMOD$.
\end{theorem}  

\vspace{6 mm}

The definition of, and the proof of statements in the theorem 
concerning $\K$, may 
be found in \cite{Bhmo} (see Theorem 3.10 there).  
Therefore  $\K$ will not appear again here.  The dual module
$\tilde{X}$ mentioned in the theorem, is discussed in 
\cite{BMP}, where we prove the analogue of the results in
pure algebra known as `Morita I' 
(see \cite{F,Bass})
That $\tilde{X} \cong Y$ will follow from Theorem 4.1
(see also 4.17 and 4.21) in \cite{BMP}, and so we
will not mention $\tilde{X}$ here again.

That strong Morita equivalence implies operator 
Morita equivalence is the easy direction of the theorem.
This follows just as in pure algebra - see \cite{BMP}
\S 3 for details.
  
One may adapt the statement of our
 main theorem above, to allow the operator
equivalence functors to be defined on not all of $OMOD$,
but only on a subcategory ${\bf D}$ of
$OMOD$ which contains $HMOD$ and the operator algebra itself,
and the maximal C$^*-$algebra it generates.
Our proof goes through 
verbatim (see comments in \cite{BlFTM}).

We remark that a number of functional 
analytic versions of the 
`Morita theorem' of equivalence of module categories, 
have been established in various contexts, although 
the categories  and methods used bear little relation to ours
(with the exception of \cite{Ri2}, which we 
will use in our proof).
We refer the reader to 
\cite{Ri2}, \cite{Beer}, and \cite{Gk}, 
for such results in the settings of 
W$^*-$algebras, unital C$^*-$algebras
and Banach algebras, respectively.
Recently in \cite{Ara1,Ara2}, Ara gave such a
Morita theorem for C$^*-$algebras, which again
is completely different to ours.  

\section{Some properties of equivalence functors}

Throughout this section $\A , \B , \C , \D $ are as before,
and $F : $ $_{\A}OMOD \rightarrow $ $_{\B}OMOD$
is an operator equivalence functor, with `inverse' $G$ (see Definition
\ref{defome}).  We set $Y = F(\A) \; , \; X = G(\B) \; , \;
Z = F(\C) \; , \;  W = F(\D)$.
In this paper we will silently be making much use of
the following two principles which are of great
assistance  with
operator algebras with c.a.i. but no identity.
Firstly, Cohen's
factorization theorem, which asserts that
a nondegenerate (left) Banach $\A$-module $X$
 has the property that $\A X = X$, and indeed
any $x \in X$ may be written as $a x'$ for
$a \in \A, x' \in X$.  Secondly, if
$\E$ is any C$^*-$algebra generated by an operator
algebra with c.a.i., then $\E$ is a nondegenerate
$\A$-module, or equivalently,
any c.a.i. for $\A$ is one for $\E$.
The latter fact is proved in \cite{BK}.
The following sequence of lemmas
will also be used extensively.  Their proofs are mostly
identical to the
analoguous results in \cite{BlFTM} and are omitted.
The first three are comparitively trivial.

\begin{lemma}
\label{coron}
Let  $V \in $ $ _{\A}OMOD$.
Then $v \mapsto r_v$ is a complete
isometry of $V$ into $_{\A}CB(\A,V)$.
The range of this map is the set $_{\A}CB^{ess}(\A,V)$.
If $V$ is also
a Hilbert space, then the map above is a completely isometric
isomorphism  $V \cong $ $_{\A}CB(\A,V)$.
\end{lemma}

\begin{lemma}
\label{bf}  If
$V, V' \in $ $_{\A}OMOD$ then the map $T \mapsto F(T)$ gives
a completely isometric
surjective linear isomorphism $_{\A}CB(V,V') \cong $ 
$_{\B}CB(F(V),F(V'))$.
If $V = V'$, then  this map is a completely isometric
isomorphism of algebras.  Moreover if $T \in \; _{\A}CB(V,V')$
is a complete isometry, then so is $F(T)$.
\end{lemma}

The last 
assertion of the previous lemma is discussed in the proof
of Theorem 8 in \cite{BMN}.

\begin{lemma}
\label{dis}  For any $V \in $ $_{\A}OMOD$, we have
$F(R_m(V)) \cong R_m(F(V))$ and $F(C_m(V)) \cong C_m(F(V))$
completely isometrically isomorphically, where $R_m(V)$
(resp. $C_m(V)$) is the operator module of rows (resp. 
columns) with $m$ elements from $V$ .
\end{lemma}

\begin{lemma}
\label{Hil}  The functors $F$ and $G$ restrict to a completely
isometric functorial equivalence of the subcategories
$_{\A}HMOD$ and $_{\B}HMOD$.
\end{lemma}

\begin{corollary}
\label{res}  The functors $F$ and $G$ restrict to a
completely isometric equivalence of $_{\C}HMOD$ and 
$_{\D}HMOD$.  This restricted 
equivalence is a normal *-equivalence
in the sense of Rieffel \cite{Rieffel2}, and 
$\C$ and $\D$ are Morita equivalent in the sense of
\cite{Rieffel2} Definition 8.17.
\end{corollary}

\begin{proof}  
This is essentially Proposition 5.1 in \cite{Bdil},
together with some general 
observations in \cite{Rieffel2} (see Definition 8.17 there).
 \end{proof}

\begin{lemma}
\label{push}
For any operator $\A$-module $V$, the canonical map
$\tau_V : Y \otimes V \rightarrow F(V)$ given by 
$y \otimes v \mapsto F(r_v)(y)$, is completely contractive
with respect to the Haagerup tensor norm,
and has dense range.  
\end{lemma}

\begin{proof}
To show $\tau_V$ has dense range,
we suppose the contrary, and let $Q$ be the nonzero quotient
map $F(V) \rightarrow \frac{F(V)}{ N}$, where 
$N = (Range \; \tau_V)^{\bar{}}$.
Then $G(Q) \neq 0$, so that
there exists $v \in V$
with $G(Q) w_V^{-1} r_v \neq 0$ as a map on $\A$, where
$w_V$ is the natural transformation $GF(V) \rightarrow V$.
Hence $FG(Q) F(w_V^{-1}) F(r_v) \neq 0$, and thus
$Q T F(r_{v}) \neq 0$ for some
 $T : F(V) \rightarrow F(V)$.  By Lemma \ref{bf}, $T = F(S)$
for some $S : V \rightarrow V$, so that 
$Q F(r_{v'}) \neq 0$
for $v' = S(v) \in V$.  Hence
$Q \circ \tau_V \neq 0$, which is a contradiction.

To show $\tau_V$ is contractive
it is sufficient to show that if
$\Vert [y_1 , \cdots , y_n ] \Vert < 1$ and $\Vert [v_1 , \cdots , v_n]^t
\Vert <1$,
then $
\Vert \sum_{k=1}^n F(r_{v_k})(y_k) \Vert < 1$.  Let us rewrite the last
expression.  Let $w = [v_1 , \cdots , v_n ]^t$ be regarded as a map in
$CB_{\A}(R_n(\A),V)$ via right multiplication $r_w$; then clearly
$\Vert r_w \Vert_{cb} < 1$.  By Lemma \ref{dis},
$F(R_n(\A))
\cong R_n(F(\A))$, so that we
may regard $[y_1 , \cdots , y_n]$ as an element
$u$ of $F(R_n(\A))$ of norm $< 1$.  We claim that $F(r_w)(u) =
\sum_{k=1}^n F(r_{v_k})(y_k)$.
This follows because $u = \sum_{k=1}^n F(i_k)
(y_k)$, where $i_k$ is the inclusion of $\A$ as the $k$-th entry in
$R_n(\A)$, so that
$$F(r_w)(u) = \sum_{k=1}^n F(r_w)
F(i_k)(y_k) =  \sum_{k=1}^n F(r_w i_k)(y_k)
= \sum_{k=1}^n F(r_{v_k})(y_k)   .
$$
Thus
$\Vert \sum_{k=1}^n r_{v_k}(y_k) \Vert = \Vert F(r_w)(u) \Vert
\leq \Vert F(r_w)  \Vert_{cb}  < 1$.
The complete contraction is similar.
\end{proof}  

\section{C$^*-$restrictable equivalences.}
It will be convenient to separate an `easy version'
of our main theorem.  We will say that an
operator equivalence functor
$F$ is
C$^*$-{\em restrictable}, if 
$F$ restricts to a functor from $_{\C}OMOD$ into
$_{\D}OMOD$.
In this section we prove our main theorem under 
the extra assumption that all functors concerned 
are C$^*$-restrictable.
First we attend to the easy direction of the theorem,
which now requires a little extra proof,
namely that the canonical equivalence functors
which come from a strong Morita equivalence, are
C$^*$-restrictable.
So suppose that  $\A$ and $\B$
are strongly Morita equivalent, and that
$X$ and $Y$ are the
strong Morita equivalence bimodules.  Then
we know from
\cite{BMN}
that $\C$ and $\D$ are strongly Morita
C$^*-$algebras, with $\D-\C-$strong Morita
equivalence bimodule
$Z \cong Y \otimes_{h\A} \C$.
Set
$F(V) = Y \otimes_{h\A} V$, for  $V$ a $\C$-operator
module.  However,  $Y \otimes_{h\A} V \cong
Y \otimes_{h\A} \C \otimes_{h\C} V \cong
Z \otimes_{h\C} V$.
Hence $F$ restricted to $_{\C}OMOD$ is
equivalent to
$Z \otimes_{h\C} -$, and is thus C$^*-$restrictable.  

Conversely,  suppose
that $F$ and $G$ are C$^*-$restrictable
operator equivalence functors.
Clearly $F$ and $G$ give an
operator Morita equivalence of
$_{\C}OMOD$ and $_{\D}OMOD$, when restricted to these
subcategories, and in \cite{BlFTM} we completely
characterized such equivalences.
Set $Y = F(\A), Z = F(\C),
X = G(\B)$ and $W = G(\D)$ as before.  From 
Lemma \ref{push},
with $V = \A$, it follows that $Y$ is a
right $\A$-operator module.  Similarly
$X$ is a right $\B$-module.  
From \cite{BlFTM}
we have that
$Z , W$ are strong Morita equivalence bimodules
for $\C$ and $\D$.
From \ref{bf},  the inclusions
$\A \subset \C$ and $\B \subset \D$ give completely
isometric inclusions $Y \rightarrow Z$
and $X \rightarrow W$.

In  \cite{BlFTM} it was shown that
 $F$ takes Hilbert $\C$-modules to Hilbert $\D$-modules.
For any Hilbert $\C$-module $K$, we have
the following sequence of canonical complete isometries
$$_{\A}CB(X,K) \cong \; _{\B}CB(\B,F(K)) \cong F(K)
\cong \; _{\D}CB(\D,F(K)) \cong \; _{\C}CB(W,K) , $$
using Lemmas \ref{coron} and \ref{bf}. 
If $R$ is the composition of this sequence of maps, then
$R$ is an inverse to the
restriction map $_{\C}CB(W,K) \rightarrow $$
_{\A}CB(X,K)$.   Hence by 3.8 in \cite{Bdil},
we have $W \cong \C \otimes_{h\A} X$ completely
isometrically and as
$\C-$modules, and it is easily checked
that this isometry is a right $\B-$module map.  Similarly,
  $Z \cong \D \otimes_{h\B} Y$.

For any $\A$-operator module $V$,
using  the last fact we see that: 
$$Y \otimes_{h\A} V
\subset \D \otimes_{h\B} (Y \otimes_{h\A} V)
\cong Z \otimes_{h\A} V  \; ,
$$
completely isometrically.
On the other hand we have the following
sequence of canonical
completely contractive $\B-$module maps:
$$
Y \otimes_{h\A} V \rightarrow
F(V) \rightarrow F(\C \otimes_{h\A} V)
\cong Z \otimes_{h\C} (\C \otimes_{h\A} V)
\cong Z \otimes_{h\C} V \; \; .
$$
The first map in this sequence comes from
\ref{push}, 
the second map comes from Lemma \ref{bf}, 
and the third map comes from the
main theorem in \cite{BlFTM}.
 The composition of the maps in this sequence coincides
with the composition of complete isometries in the last
sequence.  Hence the canonical map
$Y \otimes_{h\A} V \rightarrow F(V)$
is a complete isometry, and  is
thus a completely isometric isomorphism since it has dense
range.

Finally, $\A \cong GF(\A) \cong X \otimes_{h\B} Y$, and similarly
$\B \cong Y \otimes_{h\A} X$.  
The remaining assertions of the theorem we leave 
to the reader, namely some algebraic details 
such as checking that the transformations are natural).

\vspace{4 mm}

{\bf Remark.}  There is a natural equivalence
$_{\A^*}OMOD \cong OMOD_{\A}$, via taking the
`conjugate operator module'.   In view of this, it is
reasonable to define a `two-sided' operator 
Morita equivalence
of operator algebras, in which we adjust the definition
of left operator Morita equivalence
by replacing $F$ with two functors
$F_L : $ $_{\A}OMOD \rightarrow $ $_{\B}OMOD$ and
$F_R : OMOD_{\A} \rightarrow OMOD_{\B}$, and similarly
for $\B$.  Since
$_{\A^*}OMOD \cong OMOD_{\A}$, we get a functor
$\bar{F_R} : $$_{\A^*}OMOD \rightarrow $$_{\B^*}OMOD$.
Since $_{\C}OMOD$ is a subcategory of both
$_{\A^*}OMOD$ and $OMOD_{\A}$,
it is reasonable to assume that $F_L = \bar{F_R}$ on
$_{\C}OMOD$, and that $F_L$ is C$^*-$restrictable.
Indeed, $F_L = \bar{F_R}$ for the
canonical functors $F_L = Y \otimes_{h\A} -$ and
$F_R = - \otimes_{h\A} X$ coming from a strong Morita
equivalence.
This last interesting fact we leave as an exercise.
Thus `C$^*$-restrictability' 
is a natural condition to impose.

\section{Completion of the proof of the main theorem} 

Again $\A, \B, F, G, X, Y , W , Z$ 
are as in the previous section, but now
we fix $H \in $ $_{\A}HMOD$ to be 
the Hilbert space of the universal 
representation 
of $\C$, and fix $K = F(H)$.  Then $e(\C) \subset B(H)$, where
$e(\C)$ is the enveloping von 
Neumann algebra of $\C$.  By \ref{Hil}
and \ref{res}, $F$ and $G$
restrict to an equivalence of $_{\A}HMOD$ with $_{\B}HMOD$,
and restricts further to a
normal *-equivalence of $_{\C}HMOD$ with $_{\D}HMOD$. 
By 
\cite{Rieffel2} Propositions 1.1, 1.3 and 1.6, $\D$ acts faithfully
on $K$, and if we regard $\D$ as a subset of $B(K)$ then
the weak operator closure
$\D''$ of $\D$ in $B(K)$, is W*-isomorphic to $e(\D)$.  We shall indeed 
regard $\D$ henceforth as a subalgebra of $B(K)$.  

It is important in what follows to keep in mind the canonical 
right module action of $\B$ on $X$.
$x b = F(r_b)(x)$, for $x \in X , b \in \B$, where
as in section 2,
$r_b : \B \rightarrow \B$ is the map 
$c \mapsto cb$.  By \ref{push}, $X$ is an 
operator $\A-\B$-bimodule.
Similarly, $Y$ is canonically an
operator $\B-\A-$bimodule, and 
$Z$ and $W$ are, respectively, operator $\B-\C-$ and 
$\A-\D-$modules.  Using the last assertion in \ref{bf}
the inclusion $i$ of $\A$ in $\C$ induces a completely isometric
inclusion $F(i)$ of $Y$ in $Z$.  It is easy to see that 
$F(i)$ is a $\B-\A$-bimodule map.  
We  will regard $Y$ as a
$\B-\A-$submodule of $Z$, and, similarly, $X$ as an 
$\A-\B-$submodule of $W$.   

As we saw in Lemma \ref{push}, 
there is a left $\B-$module map
$Y \otimes X \rightarrow F(X)$ defined by
$y \otimes 
x \mapsto F(r_x)(y)$.  Since $F(X) = FG(\B) \cong \B$, 
we get a left $\B$-module map
$Y \otimes X \rightarrow \B$, which we shall write as
$[\cdot]$.  
 In a similar way we get a module map 
$(\cdot) : X \otimes Y \rightarrow \A$.  
In what
follows
 we may use the same notations for the `unlinearized' bilinear
maps, so for example we may 
use the symbols $[y,x]$ for $[y \otimes x]$.
These maps $[\cdot]$ and $(\cdot)$ have natural extensions,
which are denoted by the same symbols, to maps from
$Y \otimes W \rightarrow \D$ and $X \otimes Z \rightarrow \C$
respectively.  Namely,
$[y,w]$ is defined via $\tau_W$.
These maps $[\cdot]$ and $(\cdot)$
all have dense range, by Lemma \ref{push}.

\begin{lemma}
\label{xin}  The canonical maps 
$X \rightarrow $ $_{\B}CB(Y,\B)$ and
$Y \rightarrow $ $_{\A}CB(X,A)$ 
induced by $[\cdot]$ and $(\cdot)$
respectively, are complete isometries.  Similarly,  the
extended maps $W \rightarrow $ $_{\B}CB(Y,\D)$ ,
and $Z \rightarrow $ $_{\A}CB(X,\C)$ are complete isometries.
\end{lemma}

The proof of this is identical to the proof of the analoguous
result in \cite{BlFTM}.

The following maps $\Phi : Z  \rightarrow B(H,K)$, and $\Psi : W 
\rightarrow B(K,H)$ will play a central role in the remainder of the 
proof.  Namely, $\Phi(z)(\zeta) = F(r_\zeta)(z)$, and 
$\Psi(w)(\eta) = \omega_H 
G(r_\eta)(w)$, where $\omega_H : GF(H) \rightarrow
H$ is the $\A-$module map
coming from the natural transformation $GF \cong Id$.  
Here $r_\zeta : \C \rightarrow H$ and $r_\eta : \D \rightarrow
K$.  Since $\omega_H$ is an isometric 
surjection between Hilbert space it is unitary, and hence 
is also a $\C-$module map.
It is straightforward algebra to check that:
$$
\Psi(x) \Phi(z) = (x,z) \hspace{15mm} \& \hspace{15mm}
 \Phi(y) \Psi(w) = [y,w] V \; \;
  \eqno{(\dagger)} 
$$
for all $x \in X, y \in Y , z \in Z , w \in W$, and 
where $V \in B(K)$ is a unitary operator in 
$\D'$ composed of two natural transformations.  

\begin{lemma}
\label{phci}
The map $\Phi : Z \rightarrow B(H,K)$ 
(resp. $\Psi : W \rightarrow B(K,H)$) is a 
completely isometric $\B-\C-$module
map (resp. $\A-\D-$module map).   Moreover, $\Phi(z_1)^*\Phi(z_2) \in
\C'' = e(\C)$ for all $z_1, z_2 \in Z$, and $\Psi(w_1)^*\Psi(w_2) \in
\D''$ for $w_1, w_2 \in W$.
\end{lemma}
 
\begin{proof}
This is also almost identical to the analoguous result
in \cite{BlFTM}.
One first
establishes, for example, that for $T \in $ $\C'$, we have
$\Phi(y) T = F(T) \Phi(y)$, and this gives the 
2nd commutant assertions as in \cite{BlFTM}.

We shall simply give a few steps in the calculation showing that
$\Phi$ is a complete isometry; the missing steps may be found by 
comparison with \cite{BlFTM}:
$$
\begin{array}{ccl}
\Vert [\Phi(z_{ij})] \Vert & = & \sup 
\{ \Vert [\Phi(z_{ij})(\zeta_{kl})] \Vert
: [\zeta_{kl}] \in Ball(M_m(H^c)) , m \in \N \} \\
& = & \sup 
\{ \Vert [GF(r_{\zeta_{kl}}) \; G(r_{z_{ij}}) ] \Vert: [\zeta_{kl}] 
\in Ball(M_m(H^c)) , m \in \N \} \\
& = & \sup \{ \Vert [GF(r_{\zeta_{kl}}) \; G(r_{z_{ij}}) (x_{pq}) ] 
\Vert:
[\zeta_{kl}] \in Ball(M_m(H^c)) , [x_{pq}] \in Ball(M_r(X)) \} \\ 
& = & \sup \{ \Vert [(x_{pq},z_{ij})] \Vert : [x_{pq}] \in Ball(M_r(X)) ,
r \in \N \} \\
& = & \Vert [z_{ij}] \Vert
\end{array} 
$$
where we used the last part of Lemma \ref{xin}
in the last line.  Thus $\Phi$ is a complete isometry.
\end{proof}

\begin{lemma}
\label{unt}
The unitary $V$ is in the center of the multiplier algebra
of $\D$; and  $\Phi(y) \Psi(w) \in \D$
for all $y \in Y , w \in W$.
\end{lemma}

\begin{proof}
We will use the facts stated in the first part 
of the proof of the previous lemma.  By ($\dagger$) we know 
that $\A = [\Psi(X) \Phi(Y)]^{\bar{}}$.
Hence, using the second equation in ($\dagger$), we see that
$$ \A \Psi(X) V^{-1} \Phi(Y) 
= [\Psi(X) \Phi(Y)]^{\bar{}} \Psi(X)  V^{-1} \Phi(Y) 
\subset [\Psi(X)(\Phi(Y)\Psi(X)  V^{-1})\Phi(Y)]^{\bar{}}
\subset [\Psi(X) \Phi(Y)]^{\bar{}} = \A 
 $$
If $T \in \C'$ is such that $F(T) = V^{-1}$, then
$\Psi(X) V^{-1} \Phi(Y) = \Psi(X) \Phi(Y) T$.
Thus $\A \A T \subset \A$, so that 
$\A T \subset \A$.  Since $Y = Y \A$ we have
 $\Phi(Y) T \subset \Phi(Y)$.
Thus 
$$\Phi(y) \Psi(w) = V [y,w] = V V^{-1} \Phi(y) \Psi(w) = 
V \Phi(y) T \Psi(w) \in V \Phi(Y) \Psi(W) \subset \D \; \; .
$$  
for $y \in Y, w \in W$ .  Since $[\cdot]$ has dense range
in $\D$,
we see the multiplier assertion.
\end{proof}

\begin{theorem}
\label{slide}  The quantity 
$\Psi(w)^*\Psi(w)$, which is in $\D''$ by Lemma
\ref{phci},
is actually in $\D$
for all $w \in W$; and similarly 
$\Phi(z)^*\Phi(z) \in $ $\C$ for all $z \in Z$.  
\end{theorem}

\begin{proof}
We first observe that as in \cite{BlFTM}
the natural transformations
$GF(H) \cong H$ and
$FG(K) \cong K$  imply the following equation:
$$\langle \zeta \; \vert \; \zeta \rangle =
\sup \{  \langle (\sum_{k=1}^n \Phi(y_k)^* \Phi(y_k))\zeta \;
\vert \; \zeta \rangle  : [y_1, \cdots y_n]^t \in Ball(C_n(Y)) ,
n \in \N \} \; \;   
$$
for all $\zeta \in H$.
Replacing $\zeta$ by $\Psi(w) \eta$ for $w \in W, \eta \in K$
we have, as in \cite{BlFTM}, that
$$
\langle \; \Psi(w)^* \Psi(w) \eta \; \vert \; \eta \rangle =
\sup \{  \langle \; d \eta \; \vert \; \eta \; \rangle :
d \in \D \; , \; 0 \leq d \leq \Psi(d)^* \Psi(d) \} 
$$
A similar argument shows that  for $z \in Z, \zeta \in H$, we have
 $$
\langle \; \Phi(z)^* \Phi(z) \zeta \; \vert \; \zeta \rangle =
\sup \{  \langle \; c  \zeta \; \vert \; \zeta \rangle :
c \in \C \; , \; 0 \leq c \leq \Phi(z)^*\Phi(z)  \} 
$$
As in \cite{BlFTM} this implies that
$\Phi(z)^* \Phi(z)$ is a lowersemicontinuous
element in $e(\C) = \C''$ , for each $z \in Z$,
and that $\Psi(w)^* \Psi(w)$, as an element in
$\D''$, corresponds to a lowersemicontinuous element in $e(\D)$
(which we recall, is W$^*-$isomorphic to $\D''$).

The remainder of the proof in \cite{BlFTM} is the
same, merely replacing the `$x$' which appears in the
last few paragraphs there, by $w \in W$, and replacing the
element $a_0^2$ there by $e_\alpha^* e_\alpha$, where
$e_\alpha$ is a c.a.i. for $\A$.  We obtain
 $\Psi(w)^*\Psi(w) \in \D$.
Similarly
$\Phi(z)^* \Phi(z) \in \C$ for $z \in Z$.
\end{proof}

\begin{theorem}  
\label{stM}  The C$^*-$algebras $\C$ and $\D$ are 
strongly Morita equivalent.  In fact $Z$ , which we have seen 
to be a $\B-\C$-operator bimodule, is a $\D-\C$-strong
Morita equivalence bimodule.  Similarly, $W$ is a $\C-\D$-strong
Morita equivalence bimodule, and indeed $W \cong \bar{Z}$
unitarily (and as operator bimodules).
\end{theorem}

\begin{proof}
We will use some elementary theory or notation from 
C$^*-$modules as may be found in \cite{L2} for example.
It 
follows  by the polarization identity, and the previous
theorem,  
that $W$ is a RIGHT C$^*-$module over $\D$ with inner
product 
$\langle \; w_1 \; \vert \; w_2 \; \rangle_{\D} = 
\Psi(w_1)^*\Psi(w_2)$ .   The induced norm on $W$ from the
inner product coincides with the usual norm.   
Similarly $Z$ (or equivalently $\Phi(Z)$) is
a right C$^*-$module over $\C$.  
Also,  $W$ is a LEFT C$^*-$module over $\E =
[\Psi(W) \Psi(W)^*]^{\bar{}}$, indeed it is clear that  
$\E \cong {\mathbb K}_{\C}(Z)$, the so-called imprimitivity
C$^*-$algebra of the right C$^*-$module $Z$.
The inner product is obviously  
$_{\E}\langle \; w_1 \; \vert \; w_2 \; \rangle = 
\Psi(w_1) \Psi(w_2)^*$.
We will show that $\E = \C$.  Analoguous statements 
hold for $\D$ and $\Phi$, and we will assume below, without 
writing it down explicitly, that whenever a property
is established for $W$, the symmetric matching assertions
for $Z$.

Let $\Li$ be the linking C$^*$-algebra for 
the right C$^*-$module $W$ , viewed as a subalgebra
of $B(H \oplus K)$.  We let $\F =
[\Psi(W) \Phi(Y)]^{\bar{}}$.  It is easily seen, using
equation ($\dagger$) and Lemma \ref{unt},  that 
$\F$ is an operator algebra containing $\A$, and 
that the c.a.i. of $\A$ is a c.a.i. for $\F$.  We let
$\G = [\D \Phi(Y)]^{\bar{}}$, and we define $\M$ to be
the following subset of $B(H \oplus K)$:
  
$$
\left[ \begin{array}{cc}
\F & \Psi(W) \\
\G & \D  
\end{array}
\right] \; \; .   
$$  
This is a subalgebra by ($\dagger$) and Lemma \ref{unt}.
It is also easy to check that $\Li \M = \M$ and $\M \Li = \Li$.
Therefore from Theorem 4.15 of \cite{BMP} we conclude that
$\Li = \M$.  Comparing corners of these algebras yields
 $\E = \F$ and $G = \Psi(W)^*$.  Thus we see that 
$\A \subset \E$, from which it follows that
$\C \subset \E$, since $\C$ is the C$^*-$algebra
generated by $\A$ in $B(H)$.  Thus we have finally seen that
$W$ is a left $\C$-module, and that $\Psi$ is a 
left $\C$-module map.  By symmetry, $Z$ is a left $\D$-module
and $\Phi$ is a $\D$-module map, so that 
$$\Psi(W)^* = \G = [\D \Phi(Y)]^{\bar{}} \subset \Phi(Z) \; \; .$$ 
Also
$$\Psi(X\D) \Phi(Y) \subset [\Psi(X) \Phi(Z)]^{\bar{}}
\subset \C \; \; , 
$$
and so $\E = \F  \subset \C$.  Thus ${\mathbb K}_{\C}(Z)
\cong \E = \C$.
By symmetry note that $\Psi(W)^* = \Phi(Z)$, and that
$\D = [\Phi(Z)\Phi(Z)^*]^{\bar{}} = 
[\Psi(W)^*\Psi(W)]^{\bar{}}$.
Thus the conclusions of the theorem all hold.
\end{proof}

\begin{corollary}
\label{conn}
Operator equivalence  functors are 
automatically C$^*-$restrictable.
\end{corollary}

\begin{proof}   We keep to the notation used until 
now.  We will begin by showing
that $W$ is the maximal
dilation of $X$, and $Z$ is the maximal
dilation of $Y$.  We saw above that 
the set which we called $\G$,
 equals $Z$, so that $Y$ generates $Z$ as a
left operator $\D$-module.
We have  the following sequence 
of fairly obvious maps, using Lemmas \ref{coron} and 
\ref{bf} above:
$$_{\A}CB(X,H) \cong \; _{\B}CB(\B,K) \cong K
\cong \; _{\D}CB(\D,K) \rightarrow \; _{\A}CB(W,H) . $$
It is easily checked that  $\eta \in K$ corresponds under
the last two maps in the sequence to the map
$w \mapsto \Phi(w)(\eta)$, which lies in $_{\C}CB(W,H)$
since $\Phi$ is a left $\C$-module map.
Thus if $R$ is the composition of all the maps in 
this sequence,
then the range of $R$ is contained in
$_{\C}CB(W,K)$.  Moreover, $R$ is an inverse to the
restriction map $_{\C}CB(W,K) \rightarrow $$
_{\A}CB(X,K)$.   Thus $_{\C}CB(W,K) \cong$$
_{\A}CB(X,K)$.  Hence  by 3.8 in \cite{Bdil},
$W$ is the maximal 
dilation of $X$.  A similar argument works for $Z$.

Let $V \in \; _{\C}OMOD$.  By Lemmas \ref{coron}
and  \ref{bf} above,  and
3.9  
in \cite{Bdil}, 
we have 
$$F(V) \cong \;  _{\B}CB^{ess}(B,F(V)) \cong
\; _{\A}CB^{ess}(X,V) \cong \; _{\C}CB^{ess}(W,V) \cong
Z \otimes_{h\C} V, 
$$ 
as left $\B$-operator modules, where the 
last `$\cong$' is from \cite{Bna} Theorem
3.10.  Now the latter space is
a left $\D$-operator module, hence  
$F(V)$ is a $\D$-operator module, 
and (by a comment in \S 1, which is  
3.3 in \cite{Bdil}) 
the identity $F(V) \cong
Z \otimes_{h\C} V$ above,
 is also valid as $\D$-operator modules.  One may easily check
that this last identity is a natural isomorphism.
But $Z \otimes_{h\C} -$
is clearly a $\D$-module functor.
Hence $F$ is C$^*-$restrictable.
\end{proof}

\vspace{5 mm} 

Hence, by the result in the previous section,
our main theorem 
is proved.

\setcounter{section}{3}

\end{document}